# Enhanced energy harvesting using time-delayed feedback control from random rotational environment

Yanxia Zhang[a], Yanfei Jin[a,*], Yang Li[b]

[a] *Department of Mechanics, Beijing Institute of Technology, Beijing 100081, China*

[b] *State Key Laboratory of Mechanics and Control of Mechanical Structures, College of Aerospace Engineering, Nanjing University of Aeronautics and Astronautics, Nanjing 210016, China*

**Abstract** Motivated by improving performance of a bi-stable vibration energy harvester (VEH) from the viewpoint of vibration control, the time-delayed feedback control of displacement and velocity are constructively proposed into an electromechanical coupled VEH mounted on a rotational automobile tire, which is subject to colored noise and the periodic excitation. Using the improved stochastic averaging procedure based on energy-dependent frequency, the expressions of stationary probability density (SPD) and signal-to-noise ratio (SNR) are obtained analytically. Then, the efficiency of time-delayed feedback control on the stationary response and stochastic resonance (SR) for the delay-controlled VEH is explored in detail theoretically. Results show that both noise-induced SR and delay-induced SR can occur. Time delay is able to not only enhance the SR behavior but also weaken it. Furthermore, a larger negative feedback gain of displacement and a larger positive feedback gain of velocity are more beneficial for VEH. Interesting finding is that the optimal combination of time delay in maximizing the harvested performance, such as the harvest power, the output RMS voltage and the power conversion efficiency, is almost perfectly consistent with that in maximizing SNR. Compared with the uncontrolled VEH, the delay-controlled VEH can achieve certain desirable optimization in harvesting energy by choosing the appropriate combination of time delays and feedback gains.

## 1. Introduction

Vibration energy harvester (VEH) has the ability of utilizing ambient vibration as energy source to harvest electrical power relying on a piezoelectric cantilever in direct piezoelectric effect of vibration-to-electricity conversion. Most solutions currently for enhancing harvested energy mainly focus on resonance behavior, in which case the vibration frequency is coincident to the natural frequency of VEH [1-4]. Traditionally, linear VEH is the common design type for simplicity [5]. Nevertheless, the linear design suffers from a critical disadvantage in terms of frequency bandwidth, and ambient vibration sources mostly have wide bandwidth. Consequently, the character of narrow bandwidth in linear design makes the linear VEH insufficient to harvest energy from ambient vibration with a wider spectrum. A solution to deal with these difficulties is located in nonlinear bi-stable design for VEH which has been widely applied currently and can enhance the efficiency through increasing the relationship between ambient excitation and VEH response [6-10]. For instance, He et al. [10] has demonstrated that the nonlinear bi-stable VEH can harvest electric power efficiently with the effect of nonlinearity. Although nonlinearity has been proved as a better solution for improving performance of VEH, it is even insufficient for multiple applications related to electrical consumption necessities.

One promising potential technique for enhancing harvested energy can be achieved by using the control method, which can improve the interaction of ambient excitation and piezoelectric cantilever vibration [11,12]. From the viewpoint of vibration control, the time-delayed feedback control will be adopted to adjust nonlinear vibration and optimize system performance [13-17]. Time delay commonly exists in the feedback process of a controlled electromechanical system such that it may has significant impact on system response. Jin et al. [14,15] demonstrated in the study of a Duffing oscillator with delayed state feedback that the appropriate choice of time delay and feedback gain can enhance the control performance of dynamical systems. Yang et al. [16] studied the stiffness nonlinearities SD oscillator under time-delay control and found that time delay can not only enhance the control performance but also suppress the vibration. Furthermore, Yang et al. [17] also found in a novel hybrid energy harvester with time delay that the time-delayed feedback control can enhance stochastic resonance phenomenon leading to a large response and a high output power. Owing to the foregoing advantages of control performance, time-delayed feedback control has become a very popular solution in the field of VEH for improving the energy harvesting effectiveness [17-19]. Therefore, the time-delayed feedback control of displacement and velocity is considered in this paper and its efficiency is explored deeply for an electromechanical coupled bi-stable VEH mounted on a rotational automobile tire.

After mounted on a rotational automobile tire, the VEH can be excited autonomously by the periodic force caused by the gravity of the mass and also be disturbed inevitably by the random road excitation [20-22]. In this study, the random road excitation (i.e. road irregularity) is assumed as colored noise and generated by passing a white noise through a linear first-order filter [22-24].To the authors' knowledge, the stochastic dynamics of a bi-stable VEH with time-delayed feedback control driven by colored noise and a periodic excitation has received less attention.

Despite the challenges associated with strong nonlinearity and complexity of stochastic dynamical behavior induced by noise in bi-stable systems, it has been demonstrated [22,25,26] that the improved stochastic averaging procedure based on energy-dependent frequency is efficient to describe the Brown motion in two potential wells. There are three different periodic motions in the two potential wells according to the energy levels, i.e., vibrating in the right-side potential well, in the left-side potential well or jumping from one well to the other. Zhu et al. [26] firstly proposed the improved stochastic averaging method adopting the variable natural frequency and period of the system corresponding to the different energy levels. Later, the improved stochastic averaging procedure based on energy-dependent frequency got widely used in the study of nonlinear systems, covering the mono- and multi-stable systems [22,27,28].

The purpose of this paper is to study the stochastic dynamics in an electromechanical coupled bi-stable VEH, which is mounted on a rotational automobile tire and subjected to random road excitation. The contributions of time-delayed feedback control on enhancing the energy harvesting are discussed. The paper is organized as follows: In Sect.2, the mathematical model of delay-controlled electromechanical VEH, driven by colored noise and periodic excitation, is established and the equivalent uncoupled system is derived by introducing the harmonic transformation and integrating voltage equation. In Sect. 3, the joint SPD, the effective generalized potential function and the mean output power are obtained analytically by applying the improved stochastic averaging procedure. Then, the effects of time delay on the stochastic stationary response for delay-controlled VEH are mainly analyzed. Meanwhile, the Monte Carlo simulations (MCS) results are also given to verify the validity of the proposed theoretical method. In Sect. 4, the analytical expression of SNR is obtained, and then the effect of time delay on the delay-controlled SR is mainly explored theoretically. Meanwhile, the output RMS voltage and the power conversion efficiency affected by time delay are also analyzed in detail to evaluate the delay-controlled optimization induced by SR. Finally, some specific conclusions are drawn in Sect. 5.

## 2. Delay-controlled electromechanical VEH

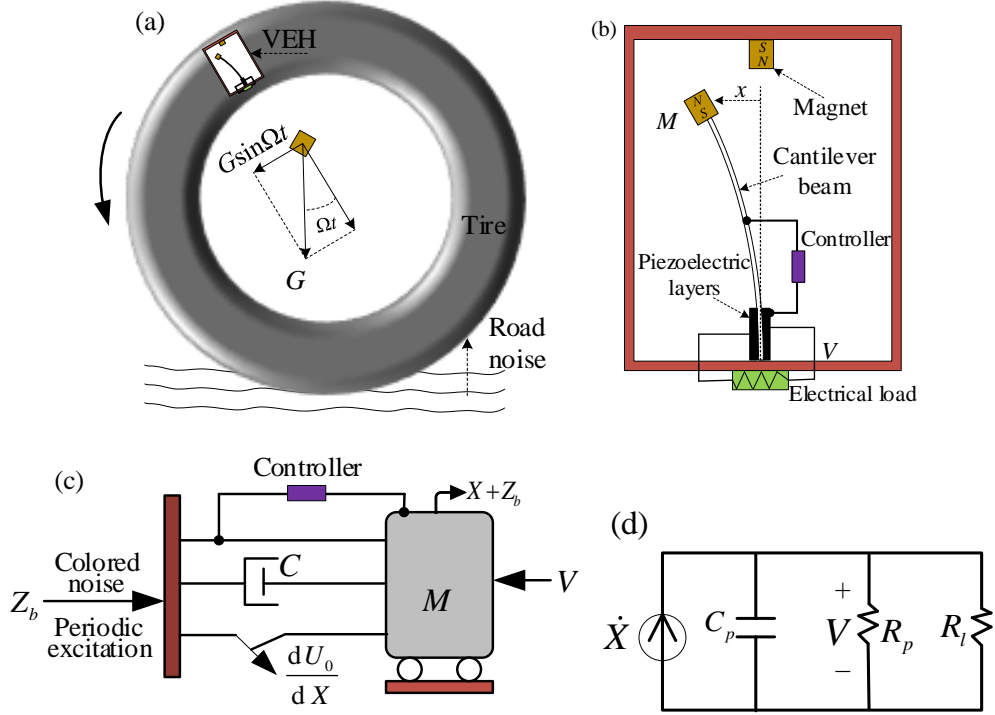

**Fig. 1** (a) A VEH mounted on a rotational automobile tire adopted from Ref. [17]; (b) Physical diagram of a coupled VEH with feedback controller; (c) Schematic of the delay-controlled VEH coupled with (d) Piezoelectric mechanism

The model of a delay-controlled electromechanical coupled VEH mounted on a rotational automobile tire is given in Fig. 1, which has been improved here by considering the time-delayed feedback control based on that of Refs. [21,22]. The VEH consists of a mechanical oscillator coupled to an electrical circuit, for vibration-to-electricity conversion by the piezoelectric mechanism, and a time-delayed feedback controller. In the case of the rotational automotive tire, the system can be autonomously driven by the periodic excitation, caused by the gravity of mass $M$, and the road irregularity which is considered as colored noise here. The controller is designed as a time-delayed feedback control of displacement and velocity. Thus, the dimensionless coupled system with time-delayed feedback control can be expressed as

$$\ddot{X} + \beta \dot{X} + \frac{dU_0(X)}{dX} + \kappa V = \mu X_{\tau_1} + \upsilon \dot{X}_{\tau_2} + \xi(t) + \varepsilon G \sin \Omega t, \qquad (1)$$
$$\dot{V} + \alpha V = \dot{X},$$

where $X_{\tau_1} = X(t-\tau_1)$, $\dot{X}_{\tau_2} = \dot{X}(t-\tau_2)$; $\tau_1$ and $\tau_2$ are the time delay; $X$ denotes the dimensionless displacement of the mass $M$; $V$ is the harvested electric voltage; $\mu$ denotes the feedback gain of dimensionless displacement; $\upsilon$ denotes the feedback gain of dimensionless velocity; $\beta$ denotes the dimensionless damping coefficient; $\kappa$ is the electromechanical coupling coefficient; $\alpha$ is the time constant ratio; $\varepsilon G \sin \Omega t$ denotes the periodic excitation caused by the gravity $G$ of mass under the rotational environment; $\varepsilon$ is the dimensionless coefficient of periodic excitation; $\xi(t)$ is the colored noise caused by road irregularity, whose statistical properties are given as

$$\langle \xi(t) \rangle = 0, \langle \xi(t)\xi(s) \rangle = \frac{D}{c}\exp\left[-\frac{|t-s|}{c}\right], \tag{2}$$

in which, $D$ and $c$ denote the noise intensity and the correlation time of colored noise, respectively.

The symmetric quartic potential function $U_0(X)$ is

$$U_0(X) = -\frac{1}{2}\delta_1 X^2 + \frac{1}{4}\delta_3 X^4, \tag{3}$$

where $\delta_1$ and $\delta_3$ denote the dimensionless linear and cubic stiffness coefficients, respectively. It has two stable equilibria and one unstable saddle point, i.e.

$$X_1^* = \sqrt{\delta_1/\delta_3},\, X_2^* = -\sqrt{\delta_1/\delta_3},\, X_3^* = 0. \tag{4}$$

By introducing the generalized harmonic function, the system displacement and the system velocity can be reformed as

$$\begin{aligned} X(t) &= A(H(t))\cos\left[\omega(H(t))t + \varphi(t)\right] + X_i^*(H(t)), \\ \dot{X}(t) &= -A(H(t))\omega(H(t))\sin\left[\omega(H(t))t + \varphi(t)\right], \end{aligned} \tag{5}$$

where $i = 1, 2, 3$ represents the three different motions vibrating in the right-side potential well, in the left-side potential well or jumping from one well to the other. $A(H)$ and $\omega(H)$ are the energy-dependent amplitude and the energy-dependent frequency, respectively. Both are slow-varying stochastic processes compared with the state variables of system.

Through the harmonic transformation, the time-delayed feedback control of displacement and velocity in Eq. (1) can be expressed approximately as

$$\begin{aligned} X(t-\tau_1) &\approx -\frac{\dot{X}(t)}{\omega(H(t))}\sin\left[\omega(H(t))\tau_1\right] + \left(X(t) - X_i^*(H(t))\right)\cos\left[\omega(H(t))\tau_1\right], \\ \dot{X}(t-\tau_2) &\approx \dot{X}(t)\cos\left[\omega(H(t))\tau_2\right] + \left(X(t) - X_i^*(H(t))\right)\omega(H(t))\sin\left[\omega(H(t))\tau_2\right]. \end{aligned} \tag{6}$$

For the purpose of uncoupling Eq. (1), the dimensionless voltage equation can be integrated and then transformed by $s = t - u$ into the following form

$$V(t) = C_1 e^{-\alpha t} + \int_0^t e^{-\alpha(t-u)}\dot{X}(u)\,\mathrm{d}u \approx \int_0^t e^{-\alpha s}\dot{X}(t-s)\,ds, \tag{7}$$

where $C_1 e^{-\alpha t}$ can be neglected for reason that the longtime stationary response is just concerned in the system.

After the harmonic transformation, similar to the derivation of Eq. (6), $\dot{X}(t-s)$ can be derived approximately as

$$\dot{X}(t-s) \approx \dot{X}(t)\cos\left[\omega(H(t))s\right] + \left(X(t) - X_i^*(H(t))\right)\omega(H(t))\sin\left[\omega(H(t))s\right]. \tag{8}$$

Consequently, by substituting Eq. (8) into Eq. (7) and neglecting the exponential decay term, the voltage can be obtained as the following form

$$V(t) = \frac{\omega^2(H(t))}{\alpha^2 + \omega^2(H(t))}\left(X(t) - X_i^*(H(t))\right) + \frac{\alpha}{\alpha^2 + \omega^2(H(t))}\dot{X}(t). \tag{9}$$

Thus, the equivalent uncoupled system can be obtained by substituting Eqs. (9) and (6) into Eq. (1)

$$\ddot{X} + \tilde{\beta}\dot{X} - \delta_1 X + \delta_3 X^3 + \tilde{\delta}(X - X_i^*) = \xi(t) + \varepsilon G \sin\Omega t, \tag{10}$$

where

$$\begin{aligned}\tilde{\beta} &= \beta + \frac{\kappa\alpha}{\alpha^2 + \omega^2(H)} + \frac{\mu}{\omega(H)}\sin(\omega(H)\tau_1) - \upsilon\cos(\omega(H)\tau_2),\\ \tilde{\delta} &= \frac{\kappa\omega^2(H)}{\alpha^2 + \omega^2(H)} - \mu\cos(\omega(H)\tau_1) - \upsilon\omega(H)\sin(\omega(H)\tau_2).\end{aligned} \tag{11}$$

From Eq. (10), one can get the quartic potential function and total energy function as

$$U(X,t) = -\frac{1}{2}\delta_1 X^2 + \frac{1}{4}\delta_3 X^4 + \frac{1}{2}\tilde{\delta}X^2 - X\varepsilon G\sin\Omega t, \tag{12}$$

$$H(X,\dot{X},t) = \frac{1}{2}\dot{X}^2 + U(X). \tag{13}$$

The energy-dependent frequency $\omega(H)$ of Eq. (11) can be obtained as

$$\omega(H) = \frac{2\pi}{T(H)} = 2\pi\left[\int_{x_a}^{x_b}\frac{dX}{\sqrt{2H - 2U(X)}}\right]^{-1}, \tag{14}$$

where $x_a$ and $x_b$ denote the min and the max displacement determined by $H = U(A)$ in the absence of delay feedback control and harmonic excitation. According to the given energy level $H$, there are three different periodic motions, i.e., vibrating in the right-side potential well, in the left-side potential well or jumping from one well to the other. It should be noted that $x_a = -x_b$ when the Brown particle jumps from one well to the other. Due to the existence of the unknown $\omega(H)$ in $U(X)$, the energy-dependent frequency $\omega(H)$ can be calculated by applying the iteration method. A similar calculative process can be found in Ref. [22].

Figure 2 shows the variation of the potential $U(X)$ in the absence of harmonic excitation from Eq. (12) and the potential well depth $\Delta U$ with different values of control parameters $\mu$, $\upsilon$, $\tau_1$ and $\tau_2$. Note that $\Delta U = |\min(U(X))|$. It can be seen that $U(X)$ is a double-well symmetric potential function about displacement $X$ (see Figs. 2(a) and 2(d)) and varies periodically with increase of $\tau_2$ as shown in Fig. 2(a). For the time-delayed velocity feedback control, as shown in Fig. 2(b), with the increase of the feedback gain $\upsilon$, well depth $\Delta U$ increases in some ranges of $\tau_2$ while decreases in other ranges. This phenomenon can also be found in Fig. 2(c) for the time-delayed displacement feedback control. These imply that the effects of feedback gains on potential function in the direction mainly depend on the chosen values of the time delays. For fixed $\tau_1 = 1.3$ and $\tau_2 = 0.5$, we find that the feedback gains of displacement and velocity can change not only the well depth but also the well space as depicted in Fig. 2(d).

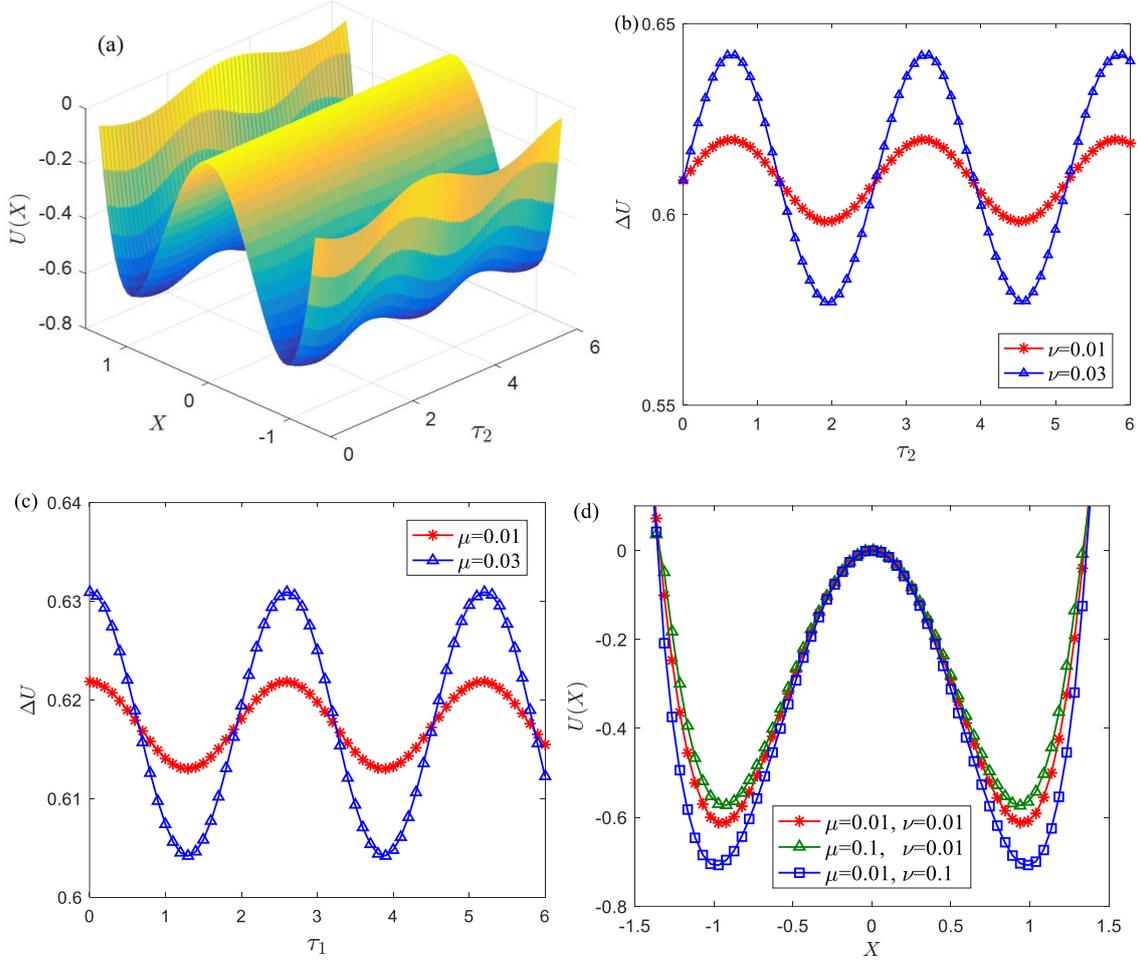

**Fig. 2** (a) The potential $U(X)$ for fixed $\mu=0.01$, $\tau_1=0.5$, $\upsilon=0.01$; (b) The function of $\Delta U$ versus $\tau_2$ with different values of $\upsilon$ for fixed $(\mu,\tau_1)=(0.01,0.5)$; (c) The function of $\Delta U$ versus $\tau_1$ with different values of $\mu$ for fixed $(\upsilon,\tau_2)=(0.01,0.5)$; (d) The variation of $U(X)$ with different values of $\mu$ and $\upsilon$ for fixed $\tau_1=1.3, \tau_2=0.5$. The other system parameters are set as $\delta_1=3$, $\delta_3=3$, $\kappa=0.3$, $\alpha=0.05$.

## 3. Delay-controlled stochastic stationary response

The equivalent uncoupled system (10) can be replaced by two first-order differential equations for $\dot{X}(t)$ and $\dot{H}(t)$ as follows

$$\begin{aligned} \dot{X} &= \pm\sqrt{2H-2U(X)}, \\ \dot{H} &= -\dot{X}\left(\tilde{\beta}\dot{X}-\tilde{\delta}X_i^*\right)+\dot{X}\xi(t). \end{aligned} \qquad (15)$$

Due to the energy process $H(t)$ being an approximate Markovian, by applying the improved stochastic averaging method of energy envelope based on the energy-dependent frequency $\omega(H)$ in Eq. (14), the averaged Itô equation for $H(t)$ can be obtained as below

$$dH = m(H)dt + \sigma(H)dB(t), \qquad (16)$$

where $B(t)$ denotes the standard Brown motion, $m(H)$ and $\sigma(H)$ denote the averaged drift and diffusion coefficients as the follows

$$m(H) = -\left\langle \dot{X}\left(\tilde{\beta}\dot{X} - \tilde{\delta}X_i^*\right)\right\rangle_t + \frac{D}{1+c^2\omega^2(H)},$$
$$\sigma^2(H) = \frac{2D}{1+c^2\omega^2(H)}\left\langle \dot{X}^2\right\rangle_t, \tag{17}$$

where $\langle \cdot \rangle_t$ denotes the time averaging, i.e.,

$$\langle \cdot \rangle_t = \frac{1}{T}\int_0^T (\cdot)\,dt = \frac{1}{T}\oint \frac{(\cdot)}{\dot{X}}\,dX = \frac{\omega(H)}{2\pi}\oint \frac{(\cdot)}{\dot{X}}\,dX. \tag{18}$$

Note that

$$\left\langle \dot{X}X_i^*\right\rangle_t = \frac{1}{T}\oint X_i^*\,dX = \frac{1}{T_R+T_L}\left[\oint_R X_R^*\,dX + \oint_L X_L^*\,dX\right] = 0, \tag{19}$$

due to the symmetry of the right-side and left-side potential wells in the absence of harmonic excitation.

Accordingly, we can obtain the SPD of the total energy from the Itô Eq.(16) as

$$p(H) = \frac{N_0}{\sigma^2(H)}\exp\left[\int \frac{2m(H)}{\sigma^2(H)}\,dH\right], \tag{20}$$

where $N_0$ is a normalization constant.

One can prove that

$$\frac{d}{dH}\ln\left[T(H)\left\langle \dot{X}^2\right\rangle_t\right] = \frac{1}{\left\langle \dot{X}^2\right\rangle_t}. \tag{21}$$

Thereby, the joint SPD of the equivalent uncoupled system (10) can be derived by $p(X,\dot{X}) = p(H)/T(H)$ as

$$p(X,\dot{X}) = \frac{N_0(1+c^2\omega^2(H))}{D}\exp\left[-\frac{\tilde{\beta}(1+c^2\omega^2(H))}{D}\left(\frac{1}{2}\dot{X}^2 - \frac{1}{2}\delta_1 X^2 + \frac{1}{4}\delta_3 X^4 + \frac{1}{2}\tilde{\delta}X^2 - X\varepsilon G\sin\Omega t\right)\right]. \tag{22}$$

From the SPD in Eq. (22), the effective generalized potential function of the system displacement and the velocity can be got as

$$\tilde{U}(X,\dot{X},t) = \tilde{\beta}(1+c^2\omega^2(H))\left(\frac{1}{2}\dot{X}^2 - \frac{1}{2}\delta_1 X^2 + \frac{1}{4}\delta_3 X^4 + \frac{1}{2}\tilde{\delta}X^2 - X\varepsilon G\cos\Omega t\right). \tag{23}$$

According to Eqs. (9) and (22), one can derive the mean-square voltage and then get the mean output power as the follows

$$E[V^2] = \int_{-\infty}^{\infty}\int_{-\infty}^{\infty}\left(\frac{\omega^2(H)}{\alpha^2+\omega^2(H)}[X-X_i^*(H)] + \frac{\alpha}{\alpha^2+\omega^2(H)}\dot{x}\right)^2 P(X,\dot{X})\,d\dot{X}\,dX, \tag{24}$$

$$E[P] = \kappa\alpha E[V^2]. \tag{25}$$

According to the analysis of the effects of time delay on the potential function (12) in the absence of harmonic excitation, the time-delayed feedback control may be of great significance in the performance optimization of energy harvester. Thus, the effects of time-delayed feedback control on the output power of the VEH (1) in the absence of harmonic excitation are analyzed. The main system parameters are set as $\delta_1 = 3$, $\delta_3 = 3$, $\kappa = 0.3$, $\alpha = 0.05$, $\beta = 0.02$, $D = 0.005$, $c = 0.3$.

First of all, the theoretical results by the stochastic averaging method are shown and compared with the numerical results by MCS from the original system Eq. (1) to verify the above obtained theoretical results (22) and (25), as displayed in Fig. 3. The joint SPD in Fig. 3(a) and the mean output power $E[P]$ denoted by solid line in Fig. 3(d) are the analytical results determined by Eqs. (22) and (25), respectively. The corresponding numerical results can be seen in Figs. 3(b) and 3(d) depicted with circle symbols. It is found that they are consistent very well.

Figure 3(c) shows the variation of the mean output power $E[P]$ by Eq. (25) with the feedback gain of displacement $\mu$ and the feedback gain of velocity $\upsilon$. For the small time delay $(\tau_1, \tau_2) = (0.5, 0.5)$, $E[P]$ increases as $\upsilon$ increases while decreases as $\mu$ increases, and $E[P]$ gets certain enhancement in the combination of a negative feedback gain of displacement and a positive feedback gain of velocity, which indicate that a larger negative feedback gain $\mu$ and a larger positive feedback gain $\upsilon$ are more beneficial for the output power of the system. Furthermore, a larger negative feedback gain $\mu$ and a larger positive feedback gain $\upsilon$ are chosen, e.g., $(\mu, \upsilon) = (-0.01, 0.01)$, and the comparison of controlled and uncontrolled cases is present in Fig. 3(d). Obviously, the output power under controlled case is superior to that under uncontrolled case. These indicate that the harvested power of the delay-controlled energy harvester can be enhanced by choosing the appropriate combination of time delays and feedback gains.

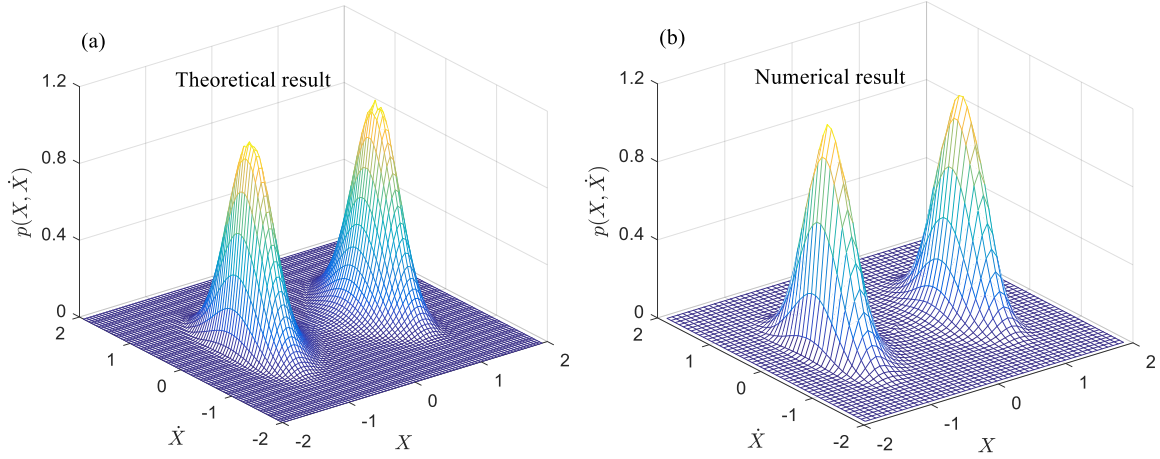

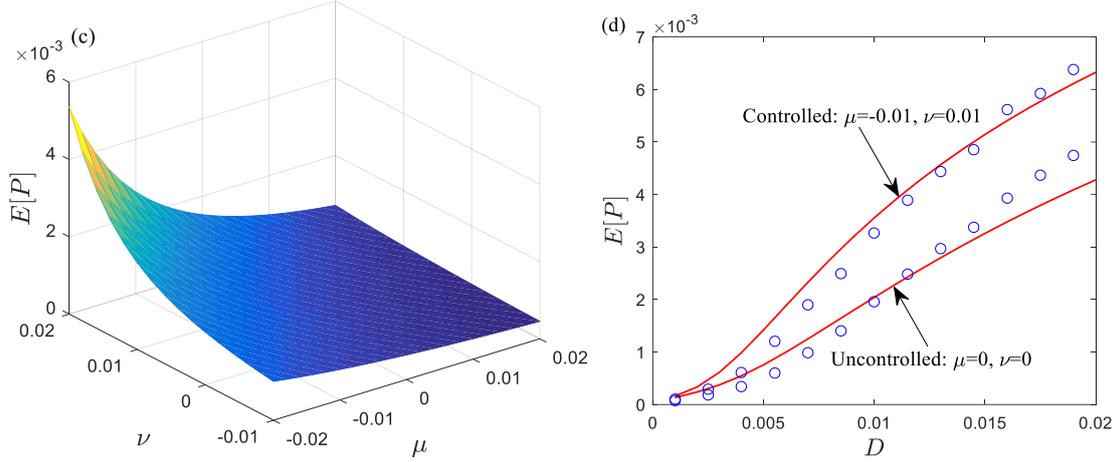

**Fig. 3** (a) The theoretical result of the joint SPD by Eq. (22) and (b) its corresponding numerical result from original system Eq. (1) for fixed $(\mu,\upsilon,\tau_1,\tau_2)=(0.01,0.01,1.3,0.5)$; (c) The dependence of the mean output power $E[P]$ on the time-delayed feedback gains $\mu$ and $\upsilon$ for fixed $(\tau_1,\tau_2)=(0.5,0.5)$; (d) The dependence of $E[P]$ on the colored noise intensity $D$ with different time-delayed feedback gains $\mu$ and $\upsilon$ for fixed $(\tau_1,\tau_2)=(0.5,0.5)$, where the solid lines denote the theoretical results and the circle symbols denote the numerical results.

The effects of time delays $\tau_1$ and $\tau_2$ on the mean output power $E[P]$ by Eq. (25) are present in Fig. 4. It can be seen in Fig. 4(a) that $E[P]$ changes periodically with the increases of $\tau_1$ and $\tau_2$. More interesting thing is that $E[P]$ has many peak values by choosing the optimal combinations of time delays $(\tau_1,\tau_2)$. For $\tau_1$ ranges in [0,2], $E[P]$ gets the peak value at $(\tau_1,\tau_2)=(0.7,2.6)$. Furthermore, the output power $E[P]$ of VEH with time-delayed feedback control of displacement and velocity at $(\mu,\upsilon,\tau_2)=(-0.005,0.005,2.6)$ outperforms that of the system only with time-delayed feedback control of displacement at $(\mu,\upsilon,\tau_2)=(-0.005,0,0)$ (see Fig. 4(b)). Meanwhile, for the controlled system with displacement feedback, a larger negative feedback gain $\mu$ of displacement can output much more power $E[P]$ markedly around $\tau_1=0.7$.

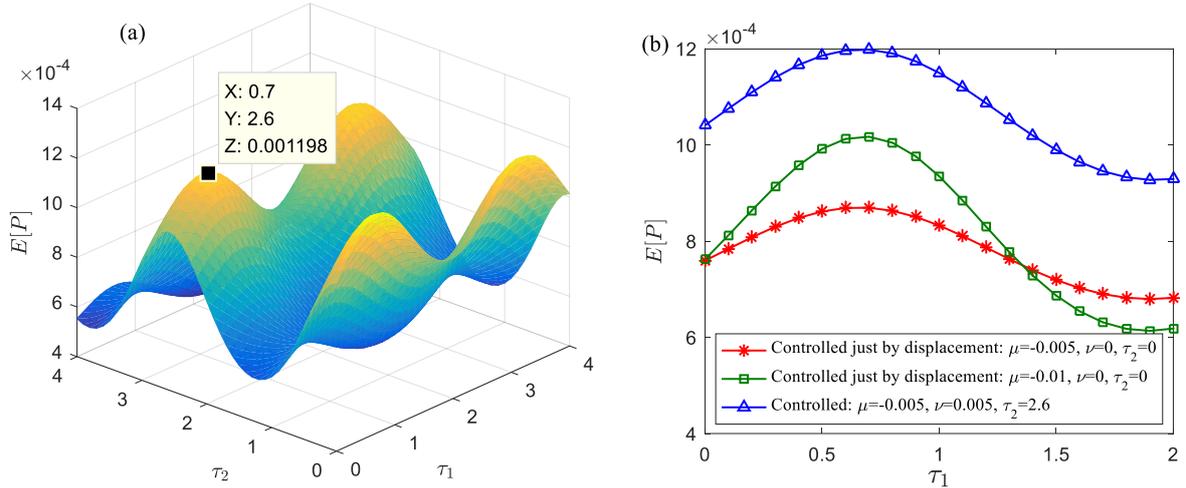

**Fig. 4** (a) The dependence of the mean output power $E[P]$ on the time delay $\tau_1$ and $\tau_2$ for fixed $(\mu,\upsilon)=(-0.005,0.005)$; (b) The dependence of $E[P]$ on the time delay $\tau_1$ with different time-delayed feedback control cases.

## 4. Delay-controlled stochastic resonance

From the above analysis, the time delay plays a constructive role in enhancing the output power in the controlled system without the harmonic excitation. For purpose of harvesting the electric energy from the external environment more effectively, it is a great idea to install the delay-controlled energy harvester on a rotational automobile tire such that it can be excited autonomously by the periodic excitation caused by the gravity of the mass. Experiments [20,21] have validated this design, showing that an energy harvester can effectively scavenge much more electric energy from rotational environment by using SR phenomenon induced by the rotational automobile tire. According to the purpose, in this section, we mainly focus on the effect of time delay on the delay-controlled VEH under the rotational environment by observing the output SNR, which can characterize the SR phenomenon quantitatively. Meanwhile, the output RMS (root mean square) voltage and the power conversion efficiency affected by time delay are also analyzed to evaluate the delay-controlled optimization capability induced by SR phenomenon. Here, the values of the periodic excitation parameters are chosen as $\varepsilon = 0.1$, $G = 0.1$, and $\Omega = 0.05$.

The SNR of the delay-controlled system can be expressed as the ratio of the output power spectrum $S_1(\Omega')$ of the signal and the output power spectrum $S_2(\Omega')$ of the noise, i.e.,

$$SNR = \frac{\int_0^\infty S_1(\Omega') d\Omega'}{S_2(\Omega' = \Omega)}. \tag{26}$$

In order to obtain the analytical expression of SNR, the Langevin equation of the equivalent uncoupled system (10) can be expressed as a two variable system

$$\dot{X} = Y,$$
$$\dot{Y} = -\tilde{\beta}\dot{X} + \delta_1 X - \delta_3 X^3 - \tilde{\delta}(X - X_i^*) + \xi(t) + \varepsilon G \sin\Omega t. \tag{27}$$

Let $\dot{X} = 0$, $\dot{Y} = 0$, two stable equilibria $X_{s_\pm}(x_{s_\pm}, 0)$ and one unstable equilibria $X_u(x_u, 0)$ of the deterministic model for system (27) without time-delayed feedback control and harmonic excitation can be obtained as

$$x_{s_\pm} = \pm\sqrt{(\delta_1 - \frac{\kappa\omega^2(H)}{\alpha^2 + \omega^2(H)})/\delta_3}, x_u = 0. \tag{28}$$

Then, the eigenvalues of the linearization characteristic equation at the three equilibria $X_{s_+}$, $X_{s_-}$ and $X_u$ can be got as

$$\lambda_m^\pm = \frac{1}{2}\left[-(c + \frac{\kappa\alpha}{\alpha^2 + \omega^2(H)}) \pm \sqrt{(c + \frac{\kappa\alpha}{\alpha^2 + \omega^2(H)})^2 - 4(-\delta_1 + \frac{\kappa\omega^2(H)}{\alpha^2 + \omega^2(H)} + 3\delta_3 x_m^2)}\right], \tag{29}$$

where $m$ denotes $s_+$, $s_-$ and $u$ in Eq. (28), respectively.
Subsequently, by using the definition of mean first-passage time and the steepest descent method, the exact expression of the transition rate $R_\pm$ out of $X_{s_\pm}$ can be obtained as

$$R_\pm = \frac{1}{2\pi}\sqrt{\frac{\lambda_{s_\pm}^+ \lambda_{s_\pm}^- \lambda_u^+}{|\lambda_u^-|}} \exp\left\{\frac{\tilde{U}(X_{s_\pm}) - \tilde{U}(X_u)}{D}\right\}. \tag{30}$$

By substituting Eq. (23) into Eq. (30), the transition rate $R_\pm$ can be further expressed as

$$R_\pm = \frac{1}{2\pi}\sqrt{\frac{\lambda_{s_\pm}^+ \lambda_{s_\pm}^- \lambda_u^+}{|\lambda_u^-|}} \exp\left\{\frac{\tilde{\beta}(1+c^2\omega^2(H))}{D}\left(-\frac{1}{2}\delta_1 x_{s_\pm}^2 + \frac{1}{4}\delta_3 x_{s_\pm}^4 + \frac{1}{2}\tilde{\delta}x_{s_\pm}^2 - x_{s_\pm}\varepsilon G \sin\Omega t\right)\right\}. \tag{31}$$

The expansion in the $\varepsilon \sin \Omega t$ term and the preservation of the first nontrivial order in the above equation (31) can lead to the following form

$$R_{\pm} = R_0 + R_1 \varepsilon \sin \Omega t, \tag{32}$$

where

$$R_0 = \frac{1}{2\pi} \sqrt{\frac{\lambda_{s_\pm}^+ \lambda_{s_\pm}^- \lambda_u^+}{|\lambda_u^-|}} \exp\left\{ \frac{\tilde{\beta}(1+c^2\omega^2(H))}{D} \left( -\frac{1}{2}\delta_1 x_{s_\pm}^2 + \frac{1}{4}\delta_3 x_{s_\pm}^4 + \frac{1}{2}\tilde{\delta} x_{s_\pm}^2 \right) \right\}, \tag{33}$$

$$R_1 = R_0 G \frac{\tilde{\beta}(1+c^2\omega^2(H))}{D} |x_{s_\pm}|. \tag{34}$$

Therefore, the output spectrum of the delay-controlled system can be obtained as

$$S(\Omega') = \frac{\pi x_{s_\pm}^2 R_1^2 \varepsilon^2}{2(R_0^2+\Omega^2)} [\delta(\Omega'-\Omega)+\delta(\Omega'+\Omega)] + \left[1 - \frac{R_1^2 \varepsilon^2}{2(R_0^2+\Omega^2)}\right] \frac{2x_{s_\pm}^2 R_0}{R_0^2+\Omega^2} = S_1(\Omega') + S_2(\Omega'). \tag{35}$$

Consequently, the analytical expression of the SNR of the delay-controlled system can be finally obtained as below by substituting Eq. (35) into Eq. (26)

$$SNR = \frac{\pi R_1^2 \varepsilon^2}{4R_0} \left[1 - \frac{R_1^2 \varepsilon^2}{2(R_0^2+\Omega^2)}\right]^{-1}. \tag{36}$$

Moreover, in order to evaluate the delay-controlled optimization capability by SR phenomenon in the energy harvester under rotational environment, the output RMS voltage $V_{rms}$ and the power conversion efficiency $\rho\%$ affected by time delay are also analyzed in the following study. Here, the power conversion efficiency $\rho\%$ is defined to assess the total efficiency from the provided mechanical power by rotational environment to the harvested electrical power, i.e.

$$\rho\% = P_e / P_m [100\%], \tag{37}$$

where $P_e = \langle \kappa \alpha V^2 \rangle$, $P_m = \langle \dot{X}\xi(t) + \dot{X}\varepsilon G \sin \Omega t \rangle$ in which $\langle \cdot \rangle$ implies time-average and ensemble average.

The effect of time delay on the output SNR are present in Fig. 5. Obviously, for fixed time delay, with the increase of noise intensity $D$, the output SNR increases firstly to a maximum and then decreases continually. This phenomenon manifests evidently that the noise-induced SR occurs. It can be observed in Fig. 5(a) that, for fixed $\tau_1 = 0.5$, as $\mu$ increases, the peak value of SNR changes slightly while the position of the peak moves toward large noise intensity gradually. Whereas, the position of the peak in SNR moves to small noise intensity gradually with the increase of $\upsilon$ for fixed $\tau_2 = 0.5$ (see Fig. 5(b)). In these two subplots, if the noise intensity $D$ of the weak noise is fixed as 0.005 displayed with the red lines, we can find that SNR increases monotonically with decreasing $\mu$ and increasing $\upsilon$. These once again indicate that a negative feedback gain $\mu$ and a positive feedback gain $\upsilon$ are more beneficial for the SR behavior of the controlled system in the case of weak noise. Whereas, in Figs. 5(c) and 5(d), for fixed $D = 0.005$ (see the red lines), the output SNR exhibits non-monotonic with the variations of time delays $\tau_1$ and $\tau_2$, which indicates that $\tau_1$ and $\tau_2$ are able to not only enhance the SR behavior but also weaken it for the delay-controlled system.

Sequentially, the joint effect of $\tau_1$ and $\tau_2$ on the output SNR is present in Fig. 5(e). We find that there exist multi peaks in the output SNR, which indicates that the SR behavior of the controlled system can be optimized by choosing the suitable combination of time delay $(\tau_1, \tau_2)$, e.g. $(\tau_1, \tau_2) = (0.6, 2.5)$. On the contrary, if chosen at the minimum of SNR, time delay can also weaken the SR behavior seriously. Interesting finding is that the optimal combination of $(\tau_1, \tau_2)$ in maximizing SNR is almost perfectly consistent with that in maximizing the mean output power $E[P]$ by comparing Figs. 5(e) and 4(a).

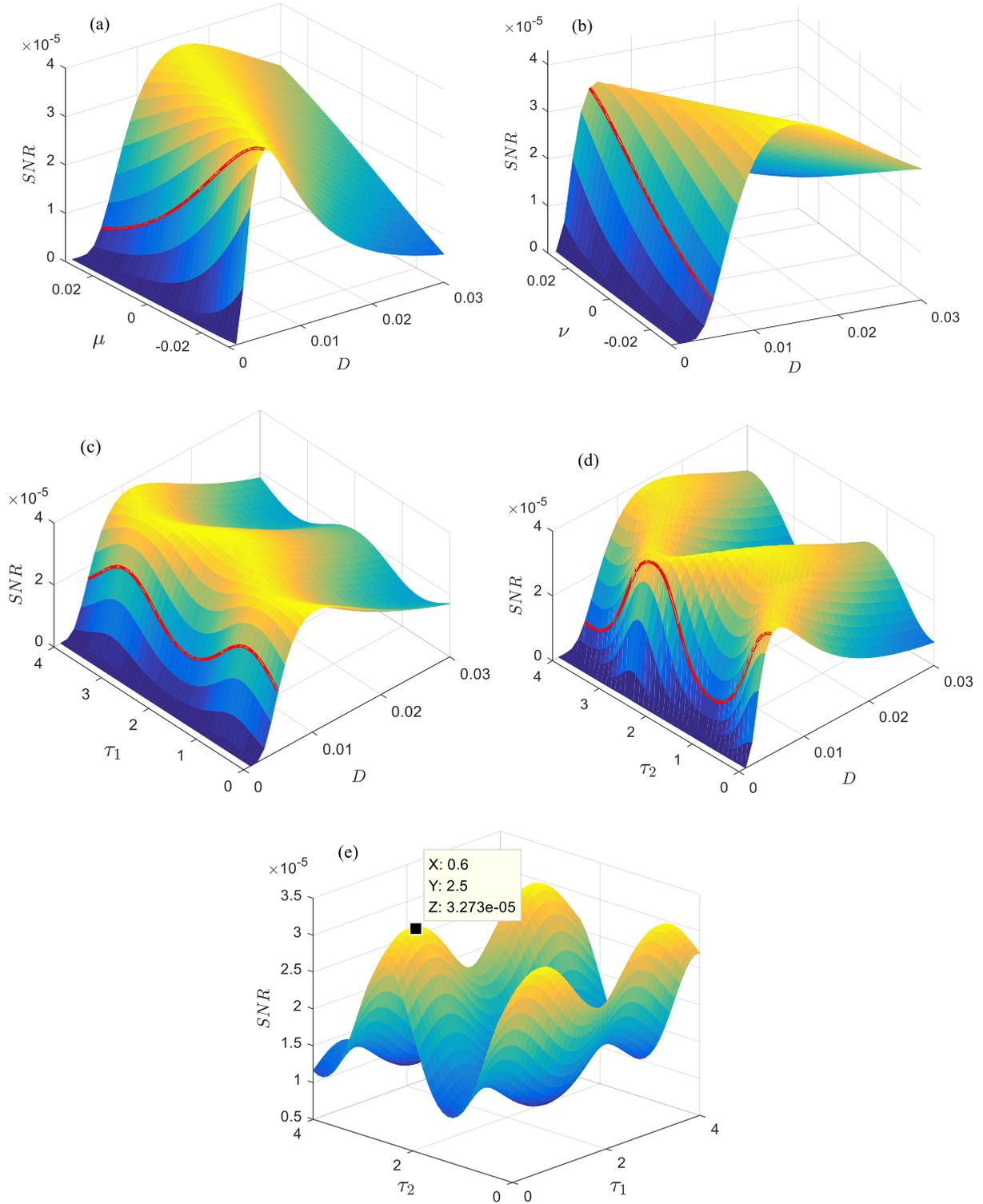

**Fig. 5**(a) SNR versus $\mu$ and $D$ for fixed $\tau_1 = 0.5$ in the absence of velocity control; (b) SNR versus $\upsilon$ and $D$ for fixed $\tau_2 = 0.5$ in the absence of displacement control; (c) SNR versus $\tau_1$ and $D$ for fixed $\mu = -0.01$ in the absence of velocity control; (d) SNR versus $\tau_2$ and $D$ for fixed $\upsilon = 0.01$ in the absence of displacement control. (e) SNR versus $\tau_1$ and $\tau_2$ for fixed $\mu = -0.005$ and $\upsilon = 0.005$. The red lines are plotted at $D = 0.005$.

Furthermore, according to the above study about the effect of time delay on the output SNR, the comparisons of delay-controlled case and uncontrolled case for the VEH in the output RMS voltage $V_{rms}$ and the power conversion efficiency $\rho\%$ are also presented in Fig. 6. For the delay-controlled case, an optimal combination of $(\tau_1,\tau_2)$ is chosen as $(0.6, 2.5)$, and the feedback gains $(\mu,\upsilon)$ are set as $(-0.01, 0.01)$. Under the time-delayed feedback control of displacement and velocity, as shown in Figs. 6(a) and 6(b), the controlled system exhibits better performance in both $V_{rms}$ and $\rho\%$ compared with the uncontrolled case, especially for the weak noise.

Figures 6(c) and 6(d) show the comparisons of controlled case and the controlled case just by displacement. As shown in Fig. 6(c), the output SNRs in both cases exhibit multi SR phenomenon as time delay $\tau_1$ increases, which indicates that delay-induced SR occurs. In addition, the curves of SNR keep the same non-monotonic movement, i.e., the positions of $\tau_1$ in maximizing SNR are the same. Meanwhile, it can be seen in Fig. 6(d) that $V_{rms}$ and $\rho\%$ reach the local maximum at the same time delay $\tau_1$ in maximizing SNR, which indicates that the system performance in the output RMS voltage and the power conversion efficiency can be optimized by delay-induced SR behavior. Moreover, for the controlled system with a positive control intensity of velocity coupled with an appropriate time delay, i.e., $(\upsilon,\tau_2)=(0.005, 2.5)$, it can be observed intuitively in both subplots that SNR, $V_{rms}$ and $\rho\%$ for the delay-controlled case are better than those for the controlled case just by displacement. Thus, the delay-controlled energy harvester can achieve certain desirable optimization in the output RMS voltage and the power conversion efficiency by choosing the appropriate combination of time delay under the rotational environment, which has a great of realistic significance in optimizing the performance for VEH.

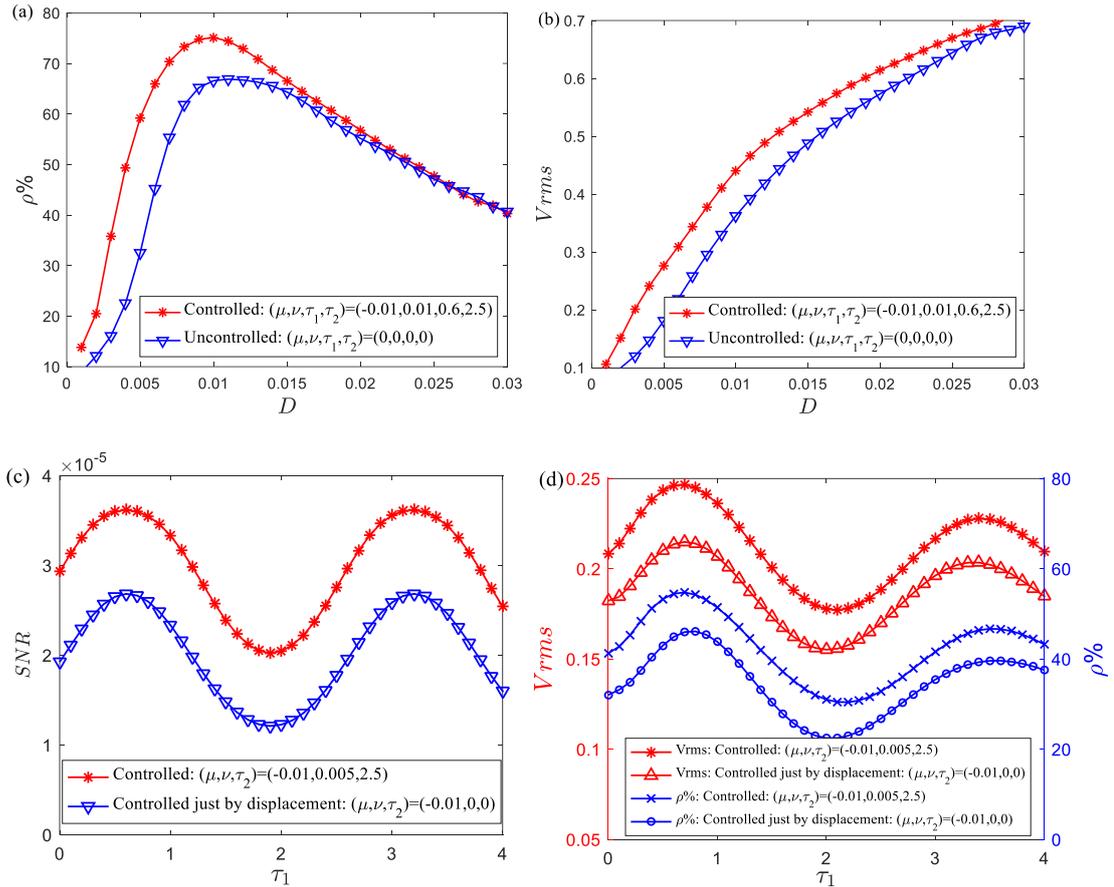

**Fig. 6** (a) The power conversion efficiency $\rho\%$ as a function of $D$ with the delay-controlled and uncontrolled cases; (b) the output RMS voltage $V_{rms}$ as a function of $D$ with the delay-controlled and uncontrolled cases; (c) the variation of SNR with the controlled case and the controlled case just by displacement; (d) the variations of $V_{rms}$ and $\rho\%$ with the controlled case and the controlled case just by displacement.

## 5. Conclusions

In this paper, for purpose of optimizing energy harvesting performance of VEH from the viewpoint of vibration control, the time-delayed feedback control of displacement and velocity are constructively proposed into an electromechanical coupled VEH mounted on a rotational automobile tire, which is driven by colored noise and periodic excitation. By applying the improved stochastic averaging procedure based on the energy-dependent frequency, the effects of time-delayed feedback control on the stochastic stationary response and SR on the bi-stable VEH are discussed based on theoretical analysis. Meanwhile, the output RMS voltage and the power conversion efficiency affected by time delays are also analyzed in detail to evaluate the delay-controlled optimization. The obtained interesting conclusions are drawn as follows:

1) The theoretical results of SPD and SNR obtained by the improved stochastic averaging procedure are well verified by the numerical results through MCS.

2) Both the mean output power $E[P]$ and the SR behavior of the delay-controlled VEH can be optimized by choosing the appropriate combination of time delays and feedback gains. A larger negative feedback gain $\mu$ of displacement and a larger positive feedback gain $\upsilon$ of velocity are more beneficial to the harvested power and the SR behavior for the controlled VEH. Time delays $\tau_1$ and $\tau_2$ can induce non-monotonic phenomenon and multi-peak appearing in the mean output power and SNR. Interesting finding is that the optimal combination of time delays in maximizing SNR is almost perfectly consistent with that in maximizing the mean output power.

3) Both noise-induced SR and delay-induced SR occur. Time delay is able to not only enhance the SR behavior but also weaken it for the delay-controlled VEH. Moreover, the output RMS voltage $V_{rms}$ and the power conversion efficiency $\rho\%$ reach the local maximum at the same time delay in maximizing SNR such that they can be optimized by delay-induced SR behavior.

4) Compared with the uncontrolled VEH, the controlled VEH exhibits better performance in the output RMS voltage and the power conversion efficiency, especially for the weak noise. All comparisons indicate that the delay-controlled VEH can achieve certain desirable optimization in harvesting energy by choosing the appropriate combination of time delay under the rotational environment, which has a great of realistic significance in optimizing the performance for VEH.

**Acknowledgements**

This study is supported by the National Natural Science Foundation of China (Nos. 11772048, 11832005) and the China Scholarship Council (CSC No. 201906030059).

**Data Availability Statement**

The data that support the findings of this study are openly available in GitHub [29].

**References**

[1] Harb, Energy harvesting: State-of-the-art, Renewable Energy 36 (10) (2011) 2641-2654.

[2] L.M. Miller, E. Halvorsen, T. Dong, P.K. Wright, Modeling and experimental verification of low-frequency MEMS energy harvesting from ambient vibrations, J. Micromech. Microeng. 21 (4) (2011) 045029.

[3] Y. Jin, S. Xiao, Y. Zhang, Enhancement of tristable energy harvesting using stochastic resonance, J. Stat. Mech. Theory E. 2018 (12) (2018) 123211.

[4] Y. Zhang, Y. Jin, P. Xu, Stochastic resonance and bifurcations in a harmonically driven tri-stable potential with colored noise, Chaos. 29 (2) (2019) 023127.

[5] S. Roundy, P.K. Wright, J.M. Rabaey, A study of low level vibrations as a power source for wireless sensor nodes, Comput. Commun. 26 (11) (2003) 1131–1144.

[6] R.L. Harne, K.W. Wang, A review of the recent research on vibration energy harvesting via bistable systems, Smart Mater. Struct.

22 (2) (2013) 023001.

[7] A. Erturk, D.J. Inman, Broadband piezoelectric power generation on high-energy orbits of the bistable duffing oscillator with electromechanical coupling, J. Sound Vib. 330 (10) (2011) 2339–2353.

[8] M.F. Daqaq, Transduction of a bistable inductive generator driven by white and exponentially correlated Gaussian noise, J. Sound Vib. 330 (11) (2011) 2554–2564.

[9] D. Liu, Y. Wu, Y. Xu, J. Li, Stochastic response of bistable vibration energy harvesting system subject to filtered Gaussian white noise, Mech. Syst. Signal Process. 130 (2019) 201-212.

[10] Q. He, M.F. Daqaq, Influence of potential function asymmetries on the performance of nonlinear energy harvesters under white noise, J. Sound Vib. 333 (15) (2014) 3479–3489.

[11] D. da Costa Ferreira, F.R. Chavarette, N.J. Peruzzi, Optimal Linear Control Driven for Piezoelectric Non-Linear Energy Harvesting on Non-Ideal Excitation Sourced, Advanced Materials Research. 971 (2014) 1107-1112.

[12] I.L. Cassidy, J.T. Scruggs, S. Behrens, Optimization of partial-state feedback for vibratory energy harvesters subjected to broadband stochastic disturbances, Smart Mater. Struct. 20 (8) (2011) 085019.

[13] D. Wu, S. Zhu, Stochastic resonance in a bistable system with time-delayed feedback and non-Gaussian noise, Phys. Lett. A 363 (3) (2007) 202–212.

[14] Y. Jin, H. Hu, Dynamics of a Duffing Oscillator with Two Time Delays in Feedback Control Under Narrow-Band Random Excitation, J Comput. Nonlin. Dyn. 3 (2) (2008) 021205.

[15] Y. Jin, H. Hu, Principal resonance of a Duffing oscillator with delayed state feedback under narrow-band random parametric excitation, Nonlinear Dyn. 50 (1-2) (2007) 213–227.

[16] T. Yang, Q. Cao, Delay-controlled primary and stochastic resonances of the SD oscillator with stiffness nonlinearities, Mech. Syst. Signal Process. 103 (2018) 216–235.

[17] T. Yang, Q. Cao, Time delay improves beneficial performance of a novel hybrid energy harvester, Nonlinear Dyn. 96 (2) (2019) 1511–1530.

[18] K.A. Alhazza, A.H. Nayfeh, M.F. Daqaq, On utilizing delayed feedback for active-multimode vibration control of cantilever beams, J. Sound Vib. 319 (3–5) (2009) 735–752.

[19] M. Hamdi, M. Belhaq, Energy harvesting in a hybrid piezoelectric-electromagnetic harvester with time delay, Recent Trends in Applied Nonlin. Mech. Phys. (2018) 69–83.

[20] Y. Zhang, R. Zheng, K. Shimono, T. Kaizuka, K. Nakano, Effectiveness testing of a piezoelectric energy harvester for an automobile wheel using stochastic resonance, Sensors. 16 (10) (2016) 1727–1742.

[21] Y. Zhang, R. Zheng, T. Kaizuka, D. Su, K. Nakano, M.P. Cartmell, Broadband vibration energy harvesting by application of stochastic resonance from rotational environments, Eur. Phys. J. Spec. Top. 224 (14-15) (2015) 2687–2701.

[22] Y. Zhang, Y. Jin, Stochastic dynamics of a piezoelectric energy harvester with correlated colored noises from rotational environment. Nonlinear Dyn. 98 (1) (2019) 501-515.

[23] S. Narayanan, S. Senthil, Stochastic optimal active control of a 2-dof quarter car model with nonlinear passive suspension elements. J. Sound Vib. 211 (3) (1998) 495–506.

[24] Y. Zhang, Y. Jin, P. Xu, S. Xiao, Stochastic bifurcations in a nonlinear tri-stable energy harvester under colored noise, Nonlinear Dyn. 99 (2020) 879–897.

[25] M. Xu, X. Li, Stochastic Averaging for Bistable Vibration Energy Harvesting System, Int. J. Mech. Sci. 141 (2018) 206-212.

[26] W. Zhu, G. Cai, R. Hu, Stochastic analysis of dynamical system with double-well potential, Int. J. Dynamics and Control 1 (1) (2013) 12-19.

[27] D. Liu, Y. Xu, J. Li, Probabilistic response analysis of nonlinear vibration energy harvesting system driven by Gaussian colored noise, Chaos Soliton. Fract. 104 (2017) 806–812.

[28] Y. Zhang, Y. Jin, P. Xu, Dynamics of a coupled nonlinear energy harvester under colored noise and periodic excitations, Int. J. Mech. Sci. 172 (2020) 105418.

[29] Y. Zhang, https://github.com/ZhangYanxia1314/Enhanced-energy-harvesting-using-time-delayed-feedback-control-from-random-rotational-environment, GitHub, 2020.